\magnification=\magstep1
\input amstex
\UseAMSsymbols
\input pictex
\vsize=23truecm
\NoBlackBoxes
\pageno=1

   \font\rmk=cmr8      \font\ttk=cmtt8

         \font\gross=cmbx10 scaled\magstep1 

 \def\mo{\operatorname{mod}}

  \def\Hom{\operatorname{Hom}}

  \def\bdim{\operatorname{\bold{dim}}}
  
 \def\arr#1#2{\arrow <1.5mm> [0.25,0.75] from #1 to #2}
\centerline{\gross The elementary 3-Kronecker modules}
     \bigskip
\centerline{Claus Michael Ringel}
           \bigskip\bigskip
{\narrower \narrower \noindent {\bf Abstract.} The 3-Kronecker quiver has two vertices,
namely a sink and a source, and 3 arrows. A regular representation of a representation-infinite
quiver such as the 3-Kronecker quiver is said to be elementary 
provided it is non-zero and not a proper extension of two regular representations. 
Of course, any regular representation has a filtration whose factors are elementary,
thus the elementary representations may be considered as the building blocks for 
obtaining all the regular representations. We are going to determine the elementary
$3$-Kronecker modules. It turns out that all the elementary modules are combinatorially
defined.\par}
           \bigskip\bigskip
Let $k$ be an algebraically closed field and $Q = K(3)$ the $3$-Kronecker quiver
  $$
	   \hbox{\beginpicture
	   \setcoordinatesystem units <1.5cm,1cm>
	   \put{$1$} at 0 0
	   \put{$2$} at 1 0
	   \arr{0.2 0.1}{0.8 0.1}
	   \arr{0.2 0}{0.8 0}
	   \arr{0.2 -.1}{0.8 -.1}
	   \endpicture}
 $$
The {\it dimension vector} of a representation $M$ of $Q$ is the pair
$(\dim M_1,\dim M_2).$ 

We denote by $A$ the arrow space of $Q$, it is a three-dimensional vector space, thus
$\Lambda = \left[\smallmatrix
k & A \cr 0 & k\endsmallmatrix\right]$ is the path algebra of $Q$.
Note that $\Lambda$ is a finite-dimensional $k$-algebra which is 
connected, hereditary and representation-infinite. The $\Lambda$-modules will be called
{\it $3$-Kronecker modules.} Of course, choosing a basis of $A$, the 3-Kronecker modules
are just the representations of $K(3)$. 

	   \medskip 
{\bf Elementary modules.} 
In general, if $\Lambda$ is a finite-dimensional $k$-algebra,
we denote by $\mo\Lambda$ the category of all (finite-dimensional left) 
$\Lambda$-modules.
We denote by $\tau$ the Auslander-Reiten translation
in $\mo\Lambda$.

Now let $\Lambda$ be the path algebra of a finite acyclic quiver.
A $\Lambda$-module $M$ is said 
to be {\it preprojective}  provided there are only finitely
many isomorphism classes of indecomposable modules $X$ with $\Hom(X,M) \neq 0$,
or, equivalently, provided $\tau^tM = 0$ for some natural number $t$. Dually,
$M$ is said 
to be {\it preinjective} provided there are only finitely
many isomorphism classes of indecomposable modules $X$ with  $\Hom(M,X)\neq 0$,
or, equivalently, provided $\tau^{-t}M = 0$ for some natural number $t$. 
A $\Lambda$-module $M$ is said to be {\it regular} provided it has no indecomposable
direct summand which is preprojective or preinjective.

A regular $\Lambda$-module $M$ is said to be {\it elementary} provided there is no
short exact sequence $0 \to M' \to M'' \to 0$ with $M', M''$ being non-zero regular 
modules (the definition is due to Crawley-Boevey, for basic results see
Kerner and Lukas [L,KL,K]) and the appendix 1.
Of course, any regular module has a filtration whose
factors are elementary. 
If $M$ is elementary, then all the modules $\tau^tM$ with $t\in \Bbb Z$
are elementary.
	\medskip
The aim of this note is to determine the elementary 3-Kronecker modules. Let
$\alpha,\beta,\gamma$ be a basis of $A$. Let $X(\alpha,\beta,\gamma)$ and  
$Y(\alpha,\beta,\gamma)$ be the 
$\Lambda$-module defined by the following pictures:
$$
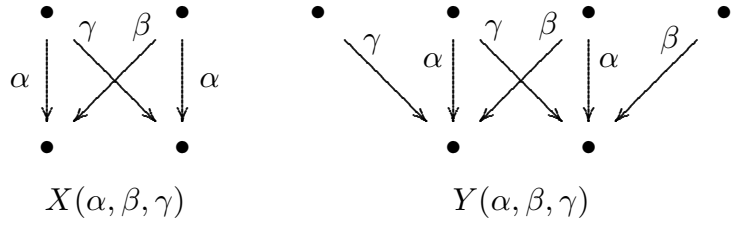
\hbox{\beginpicture
\setcoordinatesystem units <1.5cm,1.5cm>
\put{\beginpicture
 \multiput{$\bullet$} at 0 0  1 0  0 1  1 1 /
 \arr{0 0.8}{0 0.2}
 \arr{0.2 0.8}{0.8 0.2}
 \arr{1 0.8}{1 0.2}
 \arr{0.8 0.8}{0.2 0.2}
 \multiput{$\alpha$} at -.2 0.5  1.2 .5 /
 \put{$\beta$\strut} at 0.7 0.9
 \put{$\gamma$\strut} at 0.3 0.9
 \put{$X(\alpha,\beta,\gamma)$} at 0.5 -.4
 \endpicture} at 0 0 
\put{\beginpicture
 \multiput{$\bullet$} at 0 0  1 0  0 1  1 1 -1 1  2 1 /
 \arr{-.8 0.8}{-.2 0.2}
 \arr{1.8 0.8}{1.2 0.2}
 \arr{0 0.8}{0 0.2}
 \arr{0.2 0.8}{0.8 0.2}
 \arr{1 0.8}{1 0.2}
 \arr{0.8 0.8}{0.2 0.2}
 \multiput{$\alpha$} at -.15 0.65  1.15 .65 /
 \put{$\beta$\strut} at 0.7 0.9
 \put{$\gamma$\strut} at 0.3 0.9
 \put{$\beta$\strut} at 1.6 0.8
 \put{$\gamma$\strut} at -.6 0.8
 \put{$Y(\alpha,\beta,\gamma)$} at 0.5 -.4
 \endpicture} at 3 0 
\endpicture}
 $$
Here, we draw a corresponding coefficient quiver and require that all non-zero
coefficients are equal to $1$. Thus, for example $X(\alpha,\beta,\gamma) = (k^2,k^2;
\alpha,\beta,\gamma)$ with $\alpha(a,b) = (a,b)$, $\beta(a,b) = (b,0)$ and
$\gamma(a,b) = (0,a)$ for $a,b\in k.$

	\medskip
{\bf Theorem.} {\it The dimension vectors of the elementary 3-Kronecker modules are the
elements in the 
 $\tau$-orbits of $(1,1), (2,1), (2,2)$ and $(4,2)$.

Any indecomposable representation with  
dimension vector in the $\tau$-orbit of $(1,1)$ and $(2,1)$ is elementary.

An indecomposable representation with  
dimension vector $(2,2)$ or $(4,2)$ is elementary 
if and only if it is of the form $X(\alpha,\beta,\gamma)$ or
$Y(\alpha,\beta,\gamma)$, respectively for some basis $\alpha,\beta,\gamma$ of $A$.}
	\bigskip
The indecomposable representations with dimension vectors in the $\tau$-orbits of $(1,1)$
and $(2,1)$ have been studied in several papers. They are the even index Fibonacci modules,
see [FR2,FR3,R4]. 
If $M$ is indecomposable and 
$\bdim M = (1,1)$ or $(2,1)$, then there is a basis $\alpha,\beta,\gamma$ of $A$ such that
$M = B(\alpha)$ or $M = V(\beta,\gamma)$, respectively, defined as follows:
$$
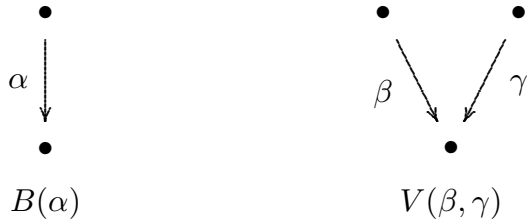
\hbox{\beginpicture
\setcoordinatesystem units <1.5cm,1.5cm>
\put{\beginpicture
 \multiput{$\bullet$} at 0 0  0 1 /
 \arr{0 0.8}{0 0.2}
 \multiput{$\alpha$} at -.2 0.5  /
 \put{$B(\alpha)$} at 0 -.4
 \endpicture} at 0 0 
\put{\beginpicture
 \multiput{$\bullet$} at 0 0  -.5 1  0.5 1   /
 \arr{-.4 0.8}{-.1 0.2}
 \arr{.4 0.8}{.1 0.2}
 \put{$\beta$\strut} at -.5 0.4
 \put{$\gamma$\strut} at 0.5 0.5
 \put{$V(\beta,\gamma)$} at 0 -.4
 \endpicture} at 3 0 
\endpicture}
 $$
Note that $B(\alpha)$ is the unique indecomposable 3-Kronecker module of length 2 which
is annihilated by $\beta$ and $\gamma$, whereas $V(\beta,\gamma)$ is the unique
indecomposable 3-Kronecker module of length 3 with simple socle which is annihilated by $\alpha$.

The indecomposable modules with dimension vector $(1,1)$ are called {\it bristles}
in [R3]. 
The indecomposable representations with dimension vector $(2,1)$ have been considered
in [BR]: there, it has been shown that any arrow $\alpha$ of a quiver gives rise to an
Auslander-Reiten sequence with indecomposable middle term say $M(\alpha)$; in this way,
we obtain the sequence:
$$
   0 \to V(\beta,\gamma) \to M(\alpha) \to \tau^{-}V(\beta,\gamma) \to 0.
$$

The study of the $\tau$-orbits of the indecomposable 3-Kronecker modules with
dimension vectors $(1,1)$ and $(2,1)$ in the papers [FR2,FR3,R4,R5] uses the universal
covering $\widetilde K(3)$ of the Kronecker quiver $K(3)$. The quiver $\widetilde K(3)$
is the 3-regular tree with bipartite orientation. Since the 3-Kronecker modules
$B(\alpha)$ and $V(\beta,\gamma)$ are cover-exceptional (they are push-downs of exceptional
representations of $\widetilde K(3)$), it follows that all the modules in the $\tau$-orbits
of $B(\alpha)$ and $V(\beta,\gamma)$ are cover-exceptional, and therefore tree modules
in the sense of [R2].

In general, one should modify the definition of a tree module as follows:
Let $Q$ be any quiver. For any pair of vertices $x,y$ of $Q$, let $A(x,y)$ be
the corresponding arrow space, this is the vector space with basis the arrows $x\to y$.
If $\alpha(1),\dots,\alpha(t)$ are the arrows $x\to y$ and 
$\beta = \sum a_i\alpha(i)$ with all $a_i\in k$ 
is an element of $A(x,y)$, we may consider for any representation
$M = (M_x,M_\alpha)_{x\in Q_0, \alpha\in Q_1}$ the linear combination
$M_\beta = \sum a_iM_{\alpha(i)}$. Given a basis $\Cal B(x,y)$ of the arrow space $A(x,y)$,
for all vertices $x,y$ of $Q$ as well as a basis
$\Cal B(M,x)$ of the vector space $M_x$, for all vertices $x$ of $Q$, 
we may write the linear maps 
$M_b$ with $b\in \Cal B(x,y)$ as matrices with respect to the bases $\Cal B(M,x),
\Cal B(M,y)$.
Looking at these matrices, we obtain a coefficient quiver $\Gamma(\Cal B(x,y),\Cal B(M,x))$
as in [R2]. A representation $M$ of the path algebra $kQ$ should be called a {\it tree module}
provided $M$ is indecomposable and there are bases $\Cal B(x,y)$ of the arrow spaces $A(x,y)$
and $\Cal B(M,x)$ of the vector spaces $M_x$ such that $\Gamma(\Cal B(x,y),\Cal B(M,x))$ is a tree.
Of course, in case $Q$ has no multiple arrows, this coincides with the definition given
in [R2]. But in general, we now allow base changes in the arrow spaces. Note that such base
changes in the arrow spaces do not effect the $kQ$-module $M$, but only its realization as
the representation of a quiver. There is the following interesting consequence: 
Any indecomposable 
representation of the 2-Kronecker quiver is a tree module, see Appendix 2.
Using this modified definition, we see immediately that {\it all the indecomposable
modules with dimension vector in the $\tau$-orbits of $(1,1)$ and $(2,1)$ are tree modules.}
On the other hand, the modules $X(\alpha,\beta,\gamma)$ (and also 
$Y(\alpha,\beta,\gamma)$) are not tree modules, see Lemma 4.2.
	\medskip 
Whereas the modules in the $\tau$-orbits of the elementary 3-Kronecker modules with
dimension vectors $(1,1)$ and $(2,1)$ are quite well understood, a
similar study of those in the $\tau$-orbits of modules with
dimension vectors $(2,2)$ and $(4,2)$ is missing. It seems that any
module $M$ in these $\tau$-orbits has a coefficient quiver with a unique cycle.
A first structure theorem for these modules is exhibited in section 5.
	\bigskip
We say that an
element $(x,y)\in K_0(\Lambda) = \Bbb Z^2$ is {\it non-negative} provided $x,y\ge 0$. 
The non-negative elements in $K_0(\Lambda)$ are just the possible dimension vectors of
$\Lambda$-modules. Note that $K_0(\Lambda)$ is endowed with the quadratic form $q$
defined by $q(x,y) = x^2+y^2-3xy$ (see for example [R1]).
A dimension vector $\bold d$ is said to be
{\it regular} provided  $q(\bold d) < 0.$ There are precisely two
$\tau$-orbits of dimension vectors $\bold d$ with $q(\bold d) = -1$, namely the 
$\tau$-orbits of $(1,1)$ and $(2,1)$. Similarly, there are precisely two
$\tau$-orbits of dimension vectors $\bold d$ with $q(\bold d) = -4$, 
namely the $\tau$-orbits of $(2,2)$ and $(4,2)$. The remaining regular dimension vectors 
$\bold d$ satisfy $q(\bold d) \le -5.$
	 \medskip 
{\bf Corollary.} {\it Let $\Lambda = kK(3)$ and $(x,y)$
a dimension vector. There exists an elementary module $M$ with dimension vector $(x,y)$ 
if and only if $q(x,y)$ is equal to $ -1$ or $-4$.}
    \bigskip
{\bf Acknowledgment.} The author wants to thank Daniel Bissinger, Rolf Farnsteiner
and Otto Kerner for fruitful discussions which initiated these investigations.
	\bigskip\bigskip 

{\bf 1\. The BGP-shift $\sigma$.}
       \medskip 
Let $\sigma$ denote the BGP-shift of $K_0(\Lambda) = \Bbb Z^2$ given
by $\sigma(x,y) = (3x-y,x),$ and let $\tau = \sigma^2.$

We denote  by $\sigma , \sigma^-$ the {\it BGP-shift functors} for
$\mo\Lambda$ (they correspond to the reflection functors of 
Bernstein-Gelfand-Ponomarev in [BGP], but take into
account that the opposite of the 3-Kronecker quiver is again the 3-Kronecker quiver).
If $M = (M_1,M_2;\alpha,\beta,\gamma)$ is a representation of $Q$, we denote by
$(\sigma  M)_1$ the kernel of the map
$\left[\smallmatrix\alpha&\beta&\gamma\endsmallmatrix\right]\:M_1^3 \to M_2$
and put $(\sigma  M)_2 = M_1$; 
the maps $\alpha,\beta,\gamma:(\sigma M)_1 \to (\sigma M)_2$ are
given by the corresponding projections. Similarly, $(\sigma^-M)_2$ is the cokernel of the map
$\left[\smallmatrix\alpha\cr\beta \cr\gamma\endsmallmatrix\right]\:M_1 \to M_2^3$ and we put
$(\sigma^-M)_1 = M_2$; now the maps $\alpha,\beta,\gamma:(\sigma^-M)_1 \to (\sigma^-M)_2$
are just the corresponding restrictions. Note that $\sigma^2$ is just the Auslander-Reiten translation $\tau$ (we should stress that this relies on the fact that we deal with a quiver
without cyclic walks of odd length, see [G]).
	\smallskip 
{\bf Remark.} The functors $\sigma$ and $\sigma^-$ 
depend
on the choice of the basis $\alpha,\beta,\gamma$ of $A$, thus we should write
$\sigma  = \sigma _{\alpha,\beta,\gamma}$ and 
$\sigma^- = \sigma^-_{\alpha,\beta,\gamma}$. 
	\medskip 
If $N$ is an indecomposable representation of $K(3)$ different from $S(2)$, then
$\bdim \sigma  N = \sigma \bdim N;$ similarly, if
$N$ is indecomposable and different from $S(1)$, then
$\bdim \sigma^- N = \sigma \bdim N$ (here, $S(1)$ and $S(2)$ are the simple
representations of $K(3)$; they are defined by
$\bdim S(1) = (1,0), \bdim S(2) = (0,1)$).

An indecomposable $\Lambda$-module $M$ is regular if and only if 
all the modules $\sigma^{n}N$ and $\sigma^{-n}N$ with $n\in \Bbb N$ are nonzero.
The restriction of $\sigma $ to the full subcategory of all regular modules
is a self-equivalence with inverse $\sigma^-$ and 
a regular module $M$ is elementary if and only if $\sigma M$ is elementary.
We say that an indecomposable representation $M$ of $K(3)$ is of 
{\it $\sigma$-type} $(x,y)$
provided $\bdim M$ belongs to the $\sigma$-orbit of $(x,y)$. 
	 \bigskip
In terms of $\sigma$, the main result can be formulated as follows: 
	\medskip 
{\bf Theorem.} {\it The elementary $kK(3)$-modules are of $\sigma$-type
$(1,1)$ and $(2,2)$. 
All the indecomposable representations of $\sigma$-type $(1,1)$ are elementary
and tree modules. 
An indecomposable representation of $\sigma$-type $(2,2)$  is either elementary
or else a tree module.}
    \medskip
The tree modules with dimension vector $(2,2)$ 
are precisely the representations of the form
$$
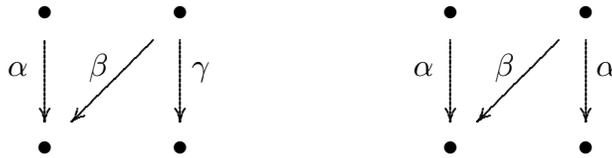
\hbox{\beginpicture
\setcoordinatesystem units <1.5cm,1.5cm>
 \put{\beginpicture
\multiput{$\bullet$} at 0 0  1 0  0 1  1 1 /
\arr{0 0.8}{0 0.2}
\arr{1 0.8}{1 0.2}
\arr{0.8 0.8}{0.2 0.2}
\put{$\alpha$\strut} at -.2 0.6
\put{$\beta$\strut} at 0.4 0.6
\put{$\gamma$\strut} at  1.15 .6
\endpicture} at 0 0
 \put{\beginpicture
\multiput{$\bullet$} at 0 0  1 0  0 1  1 1 /
\arr{0 0.8}{0 0.2}
\arr{1 0.8}{1 0.2}
\arr{0.8 0.8}{0.2 0.2}
\put{$\alpha$\strut} at -.2 0.6
\put{$\beta$\strut} at 0.4 0.6
\put{$\alpha$\strut} at  1.15 .6
\endpicture} at 3 0
\endpicture}
$$
for some basis $\alpha,\beta,\gamma$ of $A$.
           \bigskip\bigskip
 {\bf 2. Reduction to the dimension vectors $(x,y)$ with $\frac23 x \le y \le x.$}
	        \medskip
Let us denote by $\bold R$ the set of regular dimension vectors. As we have mentioned,
$\sigma$ maps $\bold R$ onto $\bold R$. There is the additional transformation $\delta$ on $K_0(\Lambda)$
defined by $\delta(x,y) = (y,x)$. Of course, it also sends $\bold R$ onto $\bold R$. If $M$ is
a representation of $Q(3)$, then $\delta(\bdim M) = \bdim M^*$, where $M^*$ is the
dual representation of $M$ (defined in the obvious way: $(M^*)_1$ is the $k$-dual of $M_2,$ $(M^*)_2$ is the $k$-dual
of $M_1$, the map $(M^*)_\alpha$ is the $k$-dual of $M_\alpha$, and so on). 
	\medskip 
{\bf Lemma.} {\it The subset 
$$
 \bold F = \{(x,y)\mid \tfrac23 x\le y \le x\}
$$
is a fundamental domain for the action of $\sigma$ and $\delta$ on $\bold R.$}
	\medskip
The proof is easy. Let us just mention that $\sigma(3,2) = (2,3)$ and that for $(x,y)\in \bold R$
with $\sigma(x,y) = (x',y')$, we have $\frac yx > \frac {y'}{x'}$ (this condition explains why we call
$\sigma$ a shift). \hfill$\square$.
	\medskip 
It follows that for looking at an elementary module, we may use the shift $\sigma$ and duality in 
order to obtain an elementary module $M$ with $\bdim M\in \bold F.$
Here is the set $\bold F$:
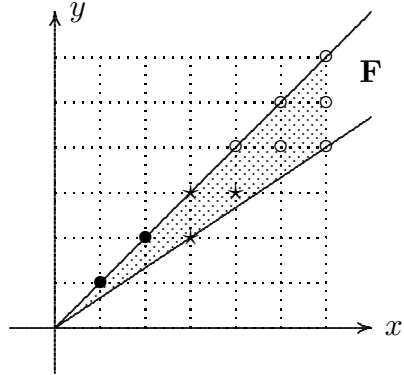
$$
 \hbox{\beginpicture
				 \setcoordinatesystem units <.5cm,.5cm>
				  \put{$\bold F$}$ at 7 5.7
				  \multiput{} at 0 -1  6 6 /
				  \arr{-1 0}{7 0}
				  \arr{0 -1}{0 7}
				  \setplotarea x from 0   to 6 , y from 0 to 6
				  \setdots <1mm>
				  \grid6 6
				  \setsolid
				  \plot 0 0  7 7 /
				  \plot 0 0  7 4.666 /
				  \multiput{$\bullet$} at 1 1  2 2 /
				  \multiput{$\star$} at 3 2  3 3  4 3 /
				  \multiput{$\circ$} at 4 4  5 4  6 4  5 5  6 5  6 6 /
				  \setshadegrid span <.5mm>
				  \vshade 0 0  0  6 4 6 /
				  \put{$x$} at 7.5 0
				  \put{$y$} at 0.5 7
 \endpicture}
 $$

 In the next section 3, we first will consider the pairs $(x,y)\in \bold F$ with $y\ge 4$,
 they are marked by a circle $\circ$. Then we deal
 with the three special pairs $(3,2), (3,3)$ and $(4,3)$ marked by a star $\star$
 (actually, instead of $(4,3)$ and $(3,2)$, we will look at $(3,4)$ and $(2,3)$, respectively).
 As we will see in section 3, all these pairs cannot occur as dimension vectors of elementary modules.

 As a consequence,
 the only possible dimension vectors in $\bold F$ which can occur as dimension vectors of elementary
 modules are $(1,1)$ and $(2,2)$; they are marked by a bullet $\bullet$ and will be studied
 in section 4.
	\bigskip\bigskip 
{\bf 3. Dimension vectors without elementary modules.} 
	\medskip
{\bf Lemma 3.1.} {\it Assume that $M$ is a regular module with a proper non-zero submodule $U$ 
such that both 
dimension vectors $\bdim U$ and $\bdim M/U$ are regular. Then $M$ is not elementary.}
         \medskip
Proof. This is a direct consequence of the fact that $M$ is elementary if and only if
for any submodule $U$ the submodule $U$ is preprojective or the factor module $M/U$ is preinjective,
see the Appendix 1. 
     \bigskip
{\bf Lemma 3.2.} {\it A $3$-Kronecker module $M$ with $\bdim M = (x,y)$ such that 
$2\le y\le x+1$ has a submodule $U$ with dimension vector $(1,2)$.}
      \medskip
Proof. Let us show that there are non-zero elements
$m\in M_1$ and $\alpha\in A$ such that $\alpha m = 0$. 
The multiplication map $A \otimes_k M_1 \to M_2$ is a linear map, let $W$ be its kernel.
Since $\dim A = 3$, we see that $\dim A\otimes_k M_1 = 3x$. Since $\dim M_2 = y,$ 
it follows that $\dim W \ge 3x-y.$ The projective
space $\Bbb P(A\otimes M_1)$ has dimension $3x-1$, the decomposable tensors in $A\otimes M_1$
form a closed subvariety $\Cal V$ of $\Bbb P(A\otimes M_1)$ of dimension $(3-1)+(x-1) = x+1$.
Since $\Cal W = \Bbb P(V)$
is a closed subspace of $\Bbb P(A\otimes M_1)$ of dimension $3x-y-1$, it follows that 
$$
 \dim(\Cal V \cap \Cal W) \ge (x+1)+(3x-y-1)-(3x-1) = x-y+1.
$$
By assumption, $x-y+1 \ge 0$, thus
$\Cal V \cap \Cal W$ is non-empty. As a consequence, we get non-zero 
elements $m\in V, \alpha\in A$ such that $\alpha m = 0,$ as required.

Given non-zero elements $m\in M_1$ and $\alpha\in A$ such that $\alpha m = 0$, 
the element $m$ generates a submodule $U'$ which is annihilated by $\alpha$, thus
$\bdim U' = (1,u)$ with $0\le u \le 2$.
Since $y\ge 2$, there is a semi-simple submodule $U''$ 
of $M$ with dimension vector $(0,2-u)$ such that
$U'\cap U'' = 0$. Let $U = U'\oplus U''$. This is a submodule of $M$ with dimension vector
$\bdim U = \bdim U'\oplus U'' = (1,2)$.
	\medskip
Remark. Under the stronger assumption $2 \le y < x$, we can argue as follows:
We have $\langle (1,2),(x,y)\rangle = x+2y - 3y = x-y > 0,$ where
$\langle-,-\rangle$ is the canonical bilinear form on $K_0(\Lambda)$ (see [R1]),
thus $\Hom(N.M) \neq 0$
for any module $N$ with $\bdim N = (1,2)$ 
The image of any non-zero map $f\:N \to M$ 
has dimension vector $(1,u)$ with $0\le u \le 2.$
       \bigskip
{\bf Lemma 3.3.} {\it If $(x,y)\in \bold F$  and $y\ge 4$, then $(x-1,y-2)$ is a regular
dimension vector.}
       \medskip
Proof. Since $y\le x,$ we have $y-2 \le x-1$.
On the other hand, the inequalities $y \ge 4$ and $y \ge \frac23 x$ imply
the inequality $y-2 \ge \frac25(x-1).$ Thus $\frac25(x-1) \le y-2 \le x-1$. As a consequence,
$(x-1,y-2)$ is a regular dimension vector. 
       \bigskip
We are now able to provide a proof for the first assertion of the Theorem:
{\it The elementary $kK(3)$-modules are of $\sigma$-type
$(1,1)$ and $(2,2)$.}
	\medskip 
Proof. Let $M$ be elementary with dimension vector $\bdim M = (x,y) \in \bold F.$
First, assume that $y \ge 4.$ According to Lemma 3.2, there is a submodule $U$ with 
the regular dimension vector $\bdim U = (1,2)$.
The factor module $M/U$ has dimension vector $(x-1,y-2)$ and $(x-1,y-1)$.  According to
Lemma 3.3, also $(x-1,y-2)$ is a regular dimension vector. Using Lemma 3.1, we obtain
a contradiction. 

It remains to show that the dimension vectors $(3,2),(3,3),(4,3)$ cannot occur.
Using duality, we may instead deal with the dimension vectors $(2,3), (3,3), (3,4)$.
Thus, assume there is given an elementary module $N$ with dimension vector   
$(2,3), (3,3)$ or $(3,4)$. According to Lemma 3.2, it has a submodule $U$ with
dimension vector $(1,2).$ The corresponding factor module $M/U$ has dimension vector
$(1,1), (2,1), (2,2)$, respectively. But all these dimension vectors are regular.
Again Lemma 3.1 provides a contradiction. 
\hfill$\square$.
	\bigskip\bigskip 
{\bf 4. The indecomposable modules with dimension vector $(1,1)$ and $(2,2)$.}
	\bigskip 
{\bf Dimension vector $(1,1)$.} {\it Any indecomposable $\Lambda$-module $M$ with dimension vector
$(1,1)$ is of the form 
$$
\hbox{\beginpicture
\setcoordinatesystem units <.8cm,.8cm>
\multiput{$\bullet$} at 0 0  0 1 /
\arr{0 0.8}{0 0.2}
\put{$\alpha$\strut} at -.3 0.6
\endpicture} 
$$
for some basis $\alpha,\beta,\gamma$ of $A$, thus a tree module.}
Namely, $M = P(1)/U$, where $P(1)$ is the indecomposable projective module corresponding to the vertex $1$
and $U$ is a two-dimensional submodule of $P(1)$. Actually, we may consider $U$ as a two-dimensional
subspace of $P(1)$. Let $\alpha,\beta,\gamma$ be a basis of $A$ such that $U = \langle \beta,\gamma\rangle$.
	\medskip 
Of course, {\it any indecomposable $\Lambda$-module with dimension vector $(1,1)$ is elementary.}
	\bigskip 
{\bf The indecomposable $\Lambda$-modules with dimension vector $(2,2)$.}
	\bigskip 
{\bf Lemma 4.1.} {\it An indecomposable module with dimension vector $(2,2)$ is elementary if and only if 
it is of the form $X(\alpha,\beta,\gamma)$.}
        \medskip
Proof. First we show: {\it The modules $M = X(\alpha,\beta,\gamma)$ are elementary.}
We have to verify that any non-zero element of $M_1$ generates a $3$-dimensional
submodule. We see this directly for the elements $(1,0)$ and $(0,1)$ of $M_1 = k^2.$ 
If $(a,b)$ with $a\neq 0, b \neq 0$, then
$\beta(a,b) = (b,0)$ and $\gamma(a,b) = (0,a)$ are linearly independent elements of $M_2 = k^2.$
This completes the proof.

Conversely, let $M$ be an elementary module with dimension vector $(2,2).$
Let us show that {\it the
restriction of $M$ to any $2$-Kronecker subalgebra has $2$-dimensional endomorphism ring.}
Let $\alpha,\beta,\gamma$ be a basis of the arrow space
and consider the restriction $M'$ of $M$ to the subquiver $K(2)$ with basis $\beta,\gamma$.
If $M'$ has a simple injective direct summand, then either $M'$ is
annihilated by $\alpha$, then $M'$ is a simple injective submodule of $M$, therefore $M$ is not
indecomposable, impossible. If $M'$ is not annihilated by $\alpha$, then $M'+\alpha(M')$ is
an indecomposable submodule of dimension $2$, thus $M$ is not elementary.
Dually, $M'$ has no simple projective direct summand. It remains to exclude the case that $M' = R\oplus R$
for some simple regular representation $R$ of $K(2)$.
Without loss of generality, we can assume that $M'$ is annihilated by $\gamma$. Since $M$ is annihilated
by $\gamma$, it is just a regular representation of the $2$-Kronecker quiver with arrow basis
$\alpha$ and $\beta$. But any 4-dimensional regular representation of a $2$-Kronecker quiver has a
$2$-dimensional regular submodule. This shows that $M$ is not elementary. Altogether we have shown that
the restriction of $M$ to any $2$-Kronecker subalgebra has $2$-dimensional endomorphism ring.

{\it If $u$ is a non-zero element of $M_1$, then $\Lambda u$ contains $M_2$ and $\dim\Lambda u = 3.$}
Namely, if $\Lambda u$ is of dimension 1, then $\Lambda u$ 
simple injective, thus $M$ cannot be indecomposable. If $\Lambda u$ is of dimension 2, then
$\Lambda u$ is a proper non-zero regular submodule and then $M$ is not elementary. It follows
that $\dim\Lambda u = 3$ and that $M_2 \subset \Lambda u.$
Given any non-zero element $u\in M_1$, there is 
a non-zero element which annihilates $u$, say $0\neq \beta\in A$. No element in 
$A\setminus\langle \beta\rangle$ annihilates $u$, since otherwise the dimension of $\Lambda u$
is at most 2.
														      Let $u, v$ be a basis of $M_1$. Let $\beta,\gamma$ be non-zero elements of $A$ with
$\beta(u) = 0, \gamma(v) = 0$. Then the elements
$\beta, \gamma$ are linearly independent, since otherwise
we would have $\gamma(u) = 0$, thus the submodule $\Lambda u$ would be of dimension at most $2$.
The elements $\beta(v), \gamma(u)$ must be linearly independent, since otherwise the restriction
of $M$ to $\beta,\gamma$ would be the direct sum of a simple projective and an indecomposable injective.
We take $\beta(v), \gamma(u)$ as an ordered basis of $M_2 = k^2,$ so that $\beta(v) = (1,0)$
and $\gamma(u) = (0,1).$
Choose an element $\alpha\in A\setminus  \langle \beta,\gamma\rangle$, thus $\alpha,\beta,\gamma$
is a basis of $A$. Let $\alpha(u) = (\kappa,\lambda)$ and $\alpha(v) = (\mu,\nu)$ 
with $\kappa,\lambda,\mu,\nu$ in $k$. 
Since $\alpha(u)$ cannot be a multiple of $\gamma(u) = (0,1),$ we see that $\kappa \neq 0$.
Since $\alpha(v)$ cannot be a multiple of $\beta(v) = (1,0),$ we see that $\nu \neq 0.$
Let $\alpha' = \alpha - \mu\beta - \lambda\gamma$.
Then
$$
\align
  \alpha'(u) &= \alpha(u)-\mu\beta(u)-\lambda\gamma(u) 
  = (\kappa,\lambda)- (0,0) - \lambda(0,1) = (\kappa,0),\cr
 \alpha'(v) &= \alpha(v)-\mu\beta(v)-\lambda\gamma(v) = (\mu,\nu)- \mu(1,0) - (0,0) =  (0,\nu).
\endalign
$$
Let $\beta' = \kappa\beta$ and $\gamma' = \nu\gamma.$ Then we have
$$
\align
   \beta'(u) = (0,0), &\quad \beta'(v) = (\kappa,0), \cr
   \gamma'(u) = (0,\nu), &\quad \gamma'(v) = (0,0),
\endalign
$$
														    Altogether, we see that
$$
\hbox{\beginpicture
\setcoordinatesystem units <1.5cm,1.5cm>
\put{$(\kappa,0)$} at 0 0
\put{$(0,\nu)$} at 1 0
\put{$u$} at  0 1
\put{$v$} at  1 1
\arr{0 0.8}{0 0.2}
\arr{0.2 0.8}{0.8 0.2}
\arr{1 0.8}{1 0.2}
\arr{0.8 0.8}{0.2 0.2}
\multiput{$\alpha'$\strut} at -.2 0.5  1.2 .5 /
\put{$\beta'$\strut} at 0.7 0.9
\put{$\gamma'$\strut} at 0.3 0.9
\endpicture}
$$
Since both elements $\kappa$ and $\nu$ are non-zero, the elements $(\kappa,0)$ and $(0,\nu)$
form a basis of $k^2,$ and $\alpha',\beta',\gamma'$ form a basis of $A$. Thus $M$ is isomorphic
to $X(\alpha',\beta',\gamma').$ This completes the proof.
	\bigskip
{\bf Lemma 4.2.} {\it A tree module with dimension vector $(2,2)$ cannot be elementary.}
	\medskip 
Proof. If $M$ is a tree module with dimension vector $(2,2)$, then the coefficient quiver has to
be of the form
$$
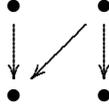
\hbox{\beginpicture
\setcoordinatesystem units <1cm,1cm>
\multiput{$\bullet$} at 0 0  0 1  1 0  1 1  /
\arr{0 0.8}{0 0.2}
\arr{1 0.8}{1 0.2}
\arr{0.8 0.8}{0.2 0.2}
\endpicture}
$$
But then $M$ has a submodule $U$ such that both $U$ and $M/U$ have dimension vector $(1,1)$.
	\bigskip
{\bf Lemma 4.3.} {\it If $M$ is indecomposable with dimension vector $(2,2)$ and not elementary, then
$M$ is of one of the following forms
$$
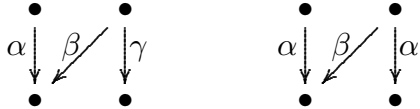
\hbox{\beginpicture
\setcoordinatesystem units <1cm,1cm>
 \put{\beginpicture
\multiput{$\bullet$} at 0 0  1 0  0 1  1 1 /
\arr{0 0.8}{0 0.2}
\arr{1 0.8}{1 0.2}
\arr{0.8 0.8}{0.2 0.2}
\put{$\alpha$\strut} at -.2 0.6
\put{$\beta$\strut} at 0.4 0.6
\put{$\gamma$\strut} at  1.15 .6
\endpicture} at 0 0
 \put{\beginpicture
\multiput{$\bullet$} at 0 0  1 0  0 1  1 1 /
\arr{0 0.8}{0 0.2}
\arr{1 0.8}{1 0.2}
\arr{0.8 0.8}{0.2 0.2}
\put{$\alpha$\strut} at -.2 0.6
\put{$\beta$\strut} at 0.4 0.6
\put{$\alpha$\strut} at  1.15 .6
\endpicture} at 3 0
\endpicture}
$$
for some basis $\alpha,\beta,\gamma$ of $A$.}
	\medskip
Proof. Let $M$ be indecomposable with dimension vector $(2,2)$. If $M$ is not
faithful, say annihilated by $0\neq \gamma \in A$, then $M$ is a $K(2)$-module
and therefore as shown on the right.

Now assume that $M$ is faithful, and not elementary.
Since $M$ is not elementary, there is an element $0\neq u \in M_1$
such that $\Lambda u$ has dimension vector $(1,1)$. The annihilator $B$ of $u$
is a 2-dimensional subspace of $A$. Let $v\in M_1\setminus \langle u\rangle$.
Since $M$ is indecomposable, we see that $\Lambda v$ has to be 3-dimensional
and there is a non-zero element $\alpha\in A$ with $\alpha v = 0.$ Since $M$
is faithful, $\alpha(u) \neq 0.$ Also, since $M$ is faithful, we have $Bv = M_2$.
Thus, there is $\beta\in B$ with $\beta(v) = \alpha(u)$. Let $\gamma \in B\setminus
\langle\beta\rangle$. Then $\alpha(u),\beta(\gamma)$ is a basis of $M_2$.
With respect to the basis $\alpha,\beta,\gamma$ of $A$, 
the basis $u,v$ of $M_1$ and the basis 
$\alpha(u),\beta(\gamma)$ of $M_2$, the module $M$ has the form as depicted on the left.
\hfill$\square$. 
        \bigskip\bigskip
{\bf 5\. The structure of the modules $\sigma^tX(\alpha,\beta,\gamma)$.}
       \medskip
The 3-Kroncker modules $I_i = \sigma^iS(2)$ are the preinjective modules, see [FR1].
	\medskip 
{\bf Proposition.} 
{\it For $t\ge 1$, there is an exact sequence}
$$
 0 \to X(\alpha,\beta,\gamma) \to \sigma^tX(\alpha,\beta,\gamma) 
    \to \bigoplus_{0\le i < t} I_i{}^2 \to 0.
$$
	\medskip
Proof. First we consider the case $t = 1.$ There is an obvious embedding of 
$X(\alpha,\beta,\gamma)$ into $Y(\alpha,\beta,\gamma) =
\sigma X(\alpha,\beta,\gamma)$, thus there is an exact sequence of the form
$$
 0 \to X(\alpha,\beta,\gamma) \to Y(\alpha,\beta,\gamma) \to S(2)^2 \to 0.
$$

Now we use induction. We start with the sequence
$$
 0 \to X(\alpha,\beta,\gamma) \to \sigma^{t}X(\alpha,\beta,\gamma) 
    \to \bigoplus_{0\le i < t} I_i{}^2 \to 0,
$$
for some $t\ge 1$ and apply $\sigma.$ In this way, we obtain the sequence
$$
 0 \to \sigma X(\alpha,\beta,\gamma) \to \sigma^{t+1}X(\alpha,\beta,\gamma) 
    \to \bigoplus_{1\le i \le t} I_i{}^2 \to 0.
$$
This shows that $M = \sigma^{t+1}X(\alpha,\beta,\gamma) $ has a submodule $U$ isomorphic
to $\sigma X(\alpha,\beta,\gamma)$, with $M/U$ isomorphic to 
$\bigoplus_{1\le i \le t} I_i{}^2$. But the case $t=1$ shows that $U$ has a submodule
$U'$ isomorphic to $X(\alpha,\beta,\gamma)$ with $U/U'$ isomorphic to $S(2)^2 = I_0{}^2$.
The embedding of $U/U'$ into  $M/U'$ has to split, since  $I_0$ is injective.
This completes the proof. 
           \bigskip\bigskip
{\bf Appendix 1. Elementary modules.}
     \medskip
According to [K], 
Proposition 4.4, a regular representation $M$ is elementary if and only if
for any
nonzero regular submodule $U$ of $M$, the factor module $M/U$ is preinjective. 
Let us include the proof of a slight improvement of this criterion. 
    \medskip 
We deal with the general setting where $\Lambda$ is a hereditary artin algebra. 
   \medskip
{\bf Proposition.} {\it Let $M$ be non-zero regular module $M$. Then $M$ is elementary 
if and only if given any submodule $U$ of $M$, 
the submodule $U$ is preprojective or the factor module $M/U$ is preinjective.}
    \medskip
Proof. Let $M$ be non-zero and regular. First, assume that for any submodule $U$ 
of $M$, the submodule $U$ is preprojective or 
the factor module $M/U$ is preinjective. Then $M$ cannot be a proper extension of 
regular modules, thus $M$ is elementary.

Conversely, let $M$ be elementary. 
Let $U$ be a submodule which is not preprojective. 
Since $M$ has no non-zero preinjective submodules, we can write $U =
U_1\oplus U_2$ with $U_1$ preprojective and $U_2$ regular. 
Since $U$ is not preprojective,
we know that $U_2$ is non-zero. Since $M$ has no non-zero preprojective  
factor modules, we decompose $M/U_2$ as a direct sum of a regular and a preinjective
module: there are submodules $V_1,V_2$ of $M$ with $V_1\cap V_2 = U, V_1+V_2 = M$
(thus $M/U = V_1/U_2\oplus V_2/U_2$) such that
$V_1/U_2$ is regular, and $V_2/U_2$ is preinjective. 

Consider $V_2$. First of
all, $V_2 \neq 0$, since $U_2$ is a non-zero submodule of $V_2$.
Second, we claim that $V_2$ is regular. Namely, $V_2$ is an
extension of the regular module $U_2$ by the preinjective module $V_2/U_2$,
thus is has no non-zero preprojective factor module. Thus, we can decompose 
$V_2 = W_1\oplus W_2$ with $W_1$ regular, $W_2$ preinjective. But $W_2$ is a
preinjective submodule of $M$, therefore $W_2 = 0.$ This shows that $V_2 = W_1$ is
regular. 

On the other hand, $W/V_2$ is isomorphic to $V_1/U_2$, thus regular. But since
$M$ is not a proper extension of regular modules, it follows that $W/V_2 = 0$,
thus $V_2 = M.$ Therefore $M/U_2 = V_2/U_2$ is preinjective. But $M/U = M/(U_1+U_2)$
is a factor module of $M/U_2$, and a factor module of a preinjective module is
preinjective. This shows that $M/U$ is preinjective. \hfill$\square$
	\bigskip
The definition of an elementary module implies that any regular module has
a filtration by elementary modules. But such filtrations are not at all unique.
This is well-known, but we would like to mention that the 3-Kronecker modules provide 
examples which are easy to remember. Here is the first such example $M$:
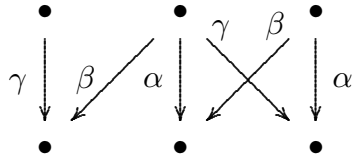
$$
\hbox{\beginpicture
\setcoordinatesystem units <1.5cm,1.5cm>
\put{\beginpicture
 \multiput{$\bullet$} at 0 0  1 0  0 1  1 1  -1 0  -1 1 /
 \arr{0 0.8}{0 0.2}
 \arr{0.2 0.8}{0.8 0.2}
 \arr{1 0.8}{1 0.2}
 \arr{0.8 0.8}{0.2 0.2}
 \arr{-.2 0.8}{-.8 0.2}
 \arr{-1 0.8}{-1 0.2}
 \multiput{$\alpha$} at -.2 0.5  1.2 .5 /
 \put{$\beta$\strut} at 0.7 0.9
 \put{$\gamma$\strut} at 0.3 0.9
 \put{$\beta$\strut} at -.7 0.5
 \put{$\gamma$\strut} at -1.2 0.5
 \endpicture} at 0 0 
\endpicture}
 $$
On the right we see that $X(\alpha,\beta,\gamma)$ is a factor module, 
and the corresponding kernel is $B(\gamma)$ (it is generated by the first 
base vector of $M_1$).
On the other hand, on the left we see that $V(\beta,\gamma)$ is a factor module,
and the corresponding kernel has dimension vector $(1,2)$
(it is generated by the last base
vector of $M_1$).

 Here is the second example $N$:
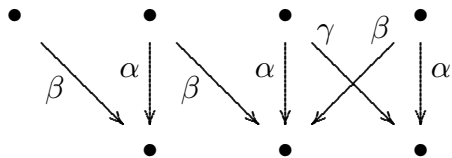
$$
\hbox{\beginpicture
\setcoordinatesystem units <1.5cm,1.5cm>
\put{\beginpicture
 \multiput{$\bullet$} at 0 0  1 0  0 1  1 1  -1 0  -1 1  -2 1 /
 \arr{0 0.8}{0 0.2}
 \arr{0.2 0.8}{0.8 0.2}
 \arr{1 0.8}{1 0.2}
 \arr{0.8 0.8}{0.2 0.2}
 \arr{-.8 0.8}{-.2 0.2}
 \arr{-1 0.8}{-1 0.2}
 \arr{-1.8 0.8}{-1.2 0.2}
 \multiput{$\alpha$} at -.15 0.6  1.15 .6  -1.15 0.6 /
 \put{$\beta$\strut} at 0.7 0.9
 \put{$\gamma$\strut} at 0.3 0.9
 \put{$\beta$\strut} at -.7 0.45
 \put{$\beta$\strut} at -1.7 0.45
 \endpicture} at 0 0 
\endpicture}
 $$
On the right we see again that $X(\alpha,\beta,\gamma)$ is a submodule, 
and the corresponding factor module is $V(\alpha,\beta)$
(generated by the first two 
base vectors of $N_1$).
On the other hand, going from left to
right, we see that the module has a filtration whose lowest two
factors are of the form 
$B(\beta)$, whereas the upper factor is $V(\alpha,\gamma)$
(generated by the last two 
base vectors of $N_1$).

    \bigskip\bigskip
\vfill\eject
{\bf Appendix 2: The representations of the $2$-Kronecker quiver.}	
	\medskip 
{\bf Proposition.} {\it Any indecomposable $K(2)$-module is a tree module (with respect to some basis of the arrow space of $K(2))$, and its coefficient quiver is of type $\Bbb A$.} 
	\medskip 
Proof. The preprojective and the preinjective modules are exceptional modules, thus they are 
tree modules with respect to any basis. 
The remaining indecomposable representations of $K(2)$ are of the form $R[t]$ 
where $R$ is simple regular, and $R[t]$ denotes the indecomposable regular
module of dimension 2t with regular socle $R$. We may choose a basis of the arrow space 
such that $R$ is isomorphic to $(k,k;1,0)$. 
Then $R[t]$ is a tree module such that the underlying graph of the coefficient quiver 
is of type $\Bbb A_{2t}$.
      \bigskip\bigskip

{\bf  References}
	\medskip 
\item{[BGP]} I.~N.~Bernstein, I.~M.~Gelfand, V.~A.~Ponomarev: Coxeter functors, and Gabriel's theorem.
		Russian mathematical surveys 28 (2) (1973), 17--32.
\item{[FR1]} Ph.~Fahr, C.~M.~Ringel: A partition formula for Fibonacci numbers.
    Journal of Integer Sequences, Vol. 11 (2008)
\item{[FR2]} Ph.~Fahr, C.~M.~Ringel: Categorification of the Fibonacci 
   Numbers Using Representations of Quivers.
   Journal of Integer Sequences. Vol. 15 (2012), Article12.2.1	
\item{[FR3]} Ph.~Fahr, C.~M.~Ringel:  The Fibonacci partition triangles.
    Advances in Mathematics 230 (2012)
\item{[G]} P.~Gabriel: Auslander-Reiten sequences and representation-finite algebras.
   Spinger LNM 831 (1980), 1--71.
\item{[K]} O.~Kerner: Representations of wild quivers.
   CMS Conf.~Proc., vol. 19. Amer.~Math. Soc., Providence, RI (1996), 65-â..107
\item{[KL]} O.~Kerner, F.~Lukas: Elementary Modules, Math.Z. 223(1996),21--434. 
\item{[L]} F.~Lukas: Elementare Moduln  \"uber wilden erblichen Algebren. Dissertation,
       D\"usseldorf (1992).
\item{[R1]} C.~M.~Ringel: Representations of $K$-species and bimodules. Representations of 
   K-species and bimodules. J. Algebra 41 (1976), 269--302.
\item{[R2]} C.~M.~Ringel: Exceptional modules are tree modules. Lin. Alg. Appl. 5-276     (1998),471--493.
\item{[R3]} C.~M.~Ringel:  Distinguished bases of exceptional modules
    In {\it Algebras, quivers and representations}. Proceedings of the Abel symposium 2011. Springer Series Abel Symposia Vol 8. (2013) 253-274. 
 
\item{[R4]} C.~M.~Ringel:
  Kronecker modules generated by modules of length 2. arXiv:1612.07679.
	\bigskip\bigskip
{\rmk
   Claus Michael Ringel \medskip
  Fakult\"at f\"ur Mathematik, Universit\"at Bielefeld\par
  D-33501 Bielefeld, Germany
    \medskip													    Department of Mathematics, Shanghai Jiao Tong University \par
  Shanghai 200240, P. R. China.
     \medskip
   e-mail: \ttk ringel\@math.uni-bielefeld.de \par}

\bye